\documentstyle[twoside]{article}
\title{The expansion of a finite number of terms of the Gauss hypergeometric function of unit argument and the Landau constants}
\author{\sc R. B. Paris\\
\\
{\em University of Abertay Dundee, Dundee DD1 1HG, UK}\\
E-Mail: r.paris@abertay.ac.uk}
\begin{document}
\def\f#1#2{\mbox{${\textstyle \frac{#1}{#2}}$}}
\def\dfrac#1#2{\displaystyle{\frac{#1}{#2}}}
\def\boldal{\mbox{\boldmath $\alpha$}}
\newcommand{\bee}{\begin{equation}}
\newcommand{\ee}{\end{equation}}
\newcommand{\lam}{\lambda}
\newcommand{\ka}{\kappa}
\newcommand{\al}{\alpha}
\newcommand{\th}{\theta}
\newcommand{\om}{\omega}
\newcommand{\Om}{\Omega}
\newcommand{\fr}{\frac{1}{2}}
\newcommand{\fs}{\f{1}{2}}
\newcommand{\g}{\Gamma}
\newcommand{\br}{\biggr}
\newcommand{\bl}{\biggl}
\newcommand{\ra}{\rightarrow}
\newcommand{\mbint}{\frac{1}{2\pi i}\int_{c-\infty i}^{c+\infty i}}
\newcommand{\mbcint}{\frac{1}{2\pi i}\int_C}
\newcommand{\mboint}{\frac{1}{2\pi i}\int_{-\infty i}^{\infty i}}
\newcommand{\gtwid}{\raisebox{-.8ex}{\mbox{$\stackrel{\textstyle >}{\sim}$}}}
\newcommand{\ltwid}{\raisebox{-.8ex}{\mbox{$\stackrel{\textstyle <}{\sim}$}}}
\renewcommand{\topfraction}{0.9}
\renewcommand{\bottomfraction}{0.9}
\renewcommand{\textfraction}{0.05}
\newcommand{\mcol}{\multicolumn}
\date{}
\maketitle
%\pagestyle{myheadings}
%\markboth{\hfill  {\it }  \hfill}
%{\hfill {\it An extension of Saalsch\"utz's theorem} \hfill}
\begin{abstract}
We obtain convergent inverse factorial expansions for the sum $S_n(a,b;c)$ of the first $n$ terms of the Gauss hypergeometric function of unit argument valid for $n\geq 1$. The form of these expansions depends on the location of the parametric excess $s:=c-a-b$ in the complex $s$-plane. The leading behaviour as $n\ra\infty$ agrees with previous results in the literature. The case $a=b=\fs$, $c=1$ corresponds to the Landau contants.

\vspace{0.4cm}

\vspace{0.3cm}

\end{abstract}

\vspace{0.3cm}

\begin{center}
{\bf 1. \  Introduction}
\end{center}
\setcounter{section}{1}
\setcounter{equation}{0}
\renewcommand{\theequation}{\arabic{section}.\arabic{equation}}
If $f(z)=\sum_{k=0}^\infty a_k z^k$ is analytic inside the unit disk and satisfies $|f(z)|<1$ when $|z|<1$, then
$\left|\sum_{k=0}^n a_k\right|\leq G_n$,
where $G_n$ are the Landau constants defined by
\bee\label{e11}
G_n=\sum_{k=0}^n 2^{-4k}\left(\!\begin{array}{c}2k\\k\end{array}\!\right)^{\!\!2}=\frac{1}{\pi}\sum_{k=0}^n\frac{\g^2(k+\fs)}{(k!)^2},\qquad n\geq 1.
\ee
It was shown by Landau \cite{L} that $G_n\sim \pi^{-1} \log\,n$ as $n\ra\infty$. 

Subsequently, it was established by Watson \cite{W} that $G_n$ is given by the convergent expansion
\bee\label{e12}
G_n=\frac{\g^2(n+\f{3}{2})}{\pi\g(n+1)\g(n+2)} \sum_{k=0}^\infty\frac{(\fs)_k (\fs)_k}{(n+2)_k k!}\{\psi(k+n+2)+\psi(k+1)-2\psi(k+\fs)\}
\ee
by writing $G_n$ as an integral over $[0,1]$ that involves the complete elliptic integral. Here and throughout
$\psi(x)=\g'(x)/\g(x)$ is the logarithmic derivative of the gamma function, also known as the $\psi$-function, and $(a)_k=\g(a+k)/\g(a)$ denotes the Pochhammer symbol. From this result combined with the basic properties $\psi(z+1)=\psi(z)+1/z$ and the asymptotic expansion \cite[p.~140]{DLMF}
\bee\label{e13a}
\psi(z)\sim \log\,z-\frac{1}{2z}-\sum_{k=1}^\infty\frac{B_{2k}}{2kz^{2k}}\qquad (z\ra+\infty)
\ee
in terms of the Bernoulli numbers $B_{2k}$, Watson deduced the asymptotic expansion 
\bee\label{e13}
G_n\sim \frac{1}{\pi} (\log\,(n+1)+\gamma+4\log\,2)-\frac{1}{4\pi(n+1)}+\frac{5}{192\pi(n+1)^2}+\cdots 
\ee
as $n\ra\infty$, where $\gamma=0.57721\ldots$ is the Euler-Mascheroni constant.

A different expansion for $G_n$ was given by Cvijovi\'c and Klinowski \cite{CK} in the form
\bee\label{e14}
G_n=\frac{1}{\pi}(\psi(n+\f{3}{2})+\gamma+4\log\,2)-\frac{1}{\pi}\sum_{k=1}^\infty \frac{(\fs)_k (\fs)_k}{k\,k! (n+\f{3}{2})_k}.
\ee
This elegant result is a convergent inverse factorial expansion that is directly amenable to computation for large $n$. Its proof is particularly simple and relies on writing the coefficients in (\ref{e11}) in terms of the Gauss hypergeometric series ${}_2F_1(\fs, \fs;k+\f{3}{2};1)$ by application of the well-known Gauss summation theorem
\bee\label{eGauss}
{}_2F_1(a,b;c;1)=\frac{\g(c) \g(c-a-b)}{\g(c-a) \g(c-b)}, \qquad \Re (c-a-b)>0
\ee
followed by expansion of the ${}_2F_1(1)$ as a convergent series. 

Based on (\ref{e14}), Nemes \cite{N} established the expansion
\bee\label{e15}
G_n\sim \frac{1}{\pi}(\log\,(n+h)+\gamma+4\log\,2)-\frac{1}{\pi}\sum_{k=1}^\infty \frac{g_k(h)}{(n+h)^k}\qquad (0<h<\f{3}{2})
\ee
as $n\ra\infty$, 
where the coefficients $g_k(h)$ are computable constants that depend on the Bernoulli polynomials and the Stirling numbers of the second kind. These coefficients are polynomials of degree $k$ in $h$ and satisfy the symmetry relation $g_k(h)=(-)^k g_k(\f{3}{2}-h)$ for $k\geq 1$. The first few $g_k(h)$ are
\[g_1(h)=\frac{1}{4}(4h-3),\quad g_2(h)=\frac{1}{192}(96h^2-144h+43),\]
\[\quad g_3(h)=\frac{1}{384}(128h^3-288h^2+172h-21).\]
Various authors have discussed the problem of determining upper and lower bounds for $G_n$; see, for example, \cite{C,CK,G,Mo,Z}.

From the second series in (\ref{e11}), it is seen that $G_n$ can be related to a finite number of terms of the series expansion of the Gauss hypergeometric function of unit argument, viz.
\[G_n={}_2F_1(\fs, \fs; 1; 1)_{\rfloor _{n+1}},\]
where the symbol $\rfloor _{n}$ signifies that only the first $n$ terms are taken. The problem of the determination
of
\[S_n(a,b;c):={}_2F_1(a,b;c;1)_{\rfloor _{n}}\]
as $n\ra\infty$ goes back to the papers of Hill \cite{H1, H2} and Bromwich \cite{Br} in the early years of the last century. By means of some lengthy algebraic manipulation and induction arguments, Hill established the leading behaviour of $S_n(a,b;c)$ for large $n$. Bromwich provided an alternative approach by use of Jensen's lemma
applied to the generalised hypergeometric series of unit argument.

In this paper, we derive convergent inverse factorial expansions for $S_n(a,b;c)$ valid for all positive integer $n$. We achieve this by means of a Mellin-Barnes integral representation, which involves routine path displacement and evaluation of residues. The type of expansion obtained depends on the location of the parameter measuring the parametric excess $s:=c-a-b$ of the hypergeometric series in the complex $s$-plane. From these results, we can obtain the asymptotic expansion (in inverse powers of $n$) of $S_n(a,b;c)$; the leading behaviour as $n\ra\infty$ is easily recovered and is found to agree with the values obtained by Hill. In the special case $a=b=\fs$, $c=1$, this leads to an alternative and simpler derivation of the inverse factorial expansion for the Landau constants $G_n$ stated in (\ref{e12}).

\vspace{0.6cm}

\begin{center}
{\bf 2. \ The expansion for $S_n(a,b;c)$}
\end{center}
\setcounter{section}{2}
\setcounter{equation}{0}
\renewcommand{\theequation}{\arabic{section}.\arabic{equation}}
We consider the sum to $n$ terms of the hypergeometric series of unit argument
\bee\label{e21}
S_n(a,b;c)=\sum_{k=0}^{n-1}\frac{(a)_k (b)_k}{(c)_k k!}
\ee
and define the associated quantities
\bee\label{e22}
s:=c-a-b,\qquad \omega_n:=\frac{\g(n+a)\g(n+b)}{\g(n)\g(n+c)},
\qquad \lambda_n:=\frac{\g(n+a)\g(n+b)}{\g(n)\g(n+a+b)}.
\ee
The parameters $a$, $b$ and $c$ are arbitrary complex constants and it will be supposed throughout that none of them equals zero or a negative integer. The quantity $s$ is known as the parametric excess; if $\Re s>0$ the series $S_n(a,b;c)$ converges to a finite limit as $n\ra\infty$ (given by Gauss' summation theorem (\ref{eGauss})), whereas if $\Re s\leq 0$ the series diverges in this limit.
From \cite[p.~81]{S}, the finite sum $S_n(a,b;c)$ can be expressed as a ${}_3F_2$ series of unit argument 
in the form
\bee\label{e23}
S_n(a,b;c)=\frac{\g(n+a)\g(n+b)}{\g(n)\g(n+a+b)}\,{}_3F_2\left[\begin{array}{c}a, b, c+n-1\\c, n+a+b\end{array};1\right].
\ee

We present the evaluation of $S_n(a,b;c)$ for real or complex values of $s$ in the following theorems.
\newtheorem{theorem}{Theorem}
\begin{theorem}$\!\!\!$.\ \ Let $n$ be a positive integer and $s=c-a-b$, with the quantities $\omega_n$ and $\lambda_n$ as defined in (\ref{e22}). Then, for finite values of $s$ such that $s\neq 0,\pm 1, \pm 2, \ldots$
and $c-a\neq 0, -1, -2, \ldots$, $c-b\neq 0, -1, -2, \ldots$ we have
\bee\label{e25}
S_n(a,b;c)=\frac{\g(c)\g(c-a-b)}{\g(c-a)\g(c-b)}-\frac{\omega_n\g(c)}{s\g(a)\g(b)}\,{}_3F_2\left[
\begin{array}{c} c-a, c-b, 1\\n+c, 1+s\end{array};1\right]. 
\ee
In the case $s=0$ we have
\[S_n(a,b;a+b)=\frac{\lambda_n\g(a+b)}{\g(a)\g(b)}\sum_{k=0}^\infty \frac{(a)_k (b)_k}{(n+a+b)_k k!}\hspace{4cm}\]
\bee\label{e26}
\hspace{2cm}\times\{\psi(n+a+b+k)+\psi(1+k)-\psi(a+k)-\psi(b+k)\}.
\ee
\end{theorem}
\vspace{0.2cm}

\noindent{\it Proof\,.}\ We employ the Mellin-Barnes integral representation for the ${}_3F_2(1)$ series
appearing in (\ref{e23}) given in \cite[p.~112]{S} in the derivation of Barnes' second lemma in the form
\[{}_3F_2\left[\begin{array}{c}a, b, c\\d, e\end{array};1\right]=\frac{\g(d)\g(e)}{\g(a)\g(b) \g(e-c)\g(d-a)\g(d-b)}\]
\[\times\frac{1}{2\pi i}\int_{-\infty i}^{\infty i}\frac{\g(a+\tau)\g(b+\tau)}{\g(e+\tau)}\,
\g(d-a-b-\tau)\g(e-c+\tau)\g(-\tau)\,d\tau,\]
where the integration path lies to the left of the poles of $\g(-\tau)$ and $\g(d-a-b-\tau)$. Then we find that
\[S_n(a,b;c)=\frac{\lambda_n \g(c)\g(n+a+b)}{\g(a)\g(b)\g(c-a)\g(c-b)\g(1-s)}\hspace{5cm}\]
\bee\label{e28}
\hspace{2cm}\times \frac{1}{2\pi i}\int_{-\infty i}^{\infty i}
\frac{\g(a+\tau)\g(b+\tau)}{\g(n+a+b+\tau)}\g(1-s+\tau)\g(s-\tau)\g(-\tau)\,d\tau,
\ee
where the integration path may be suitably indented (if necessary) to separate the poles of $\g(-\tau)$ and $\g(s-\tau)$ from those of $\g(1-s+\tau)$, $\g(a+\tau)$ and $\g(b+\tau)$. The above separation of the sequences of poles in the $\tau$-plane is possible if $s\neq 1, 2, \ldots\,$, and $c-a$, $c-b\neq 0, -1, -2, \ldots$ and $a$, $b\neq 0, -1, -2, \ldots\ $.

Provided $s$ is not an integer, the poles on the right of the path are all simple. 
We denote the integral appearing in (\ref{e28}) (including the factor $(2\pi i)^{-1}$) by $I$.
Displacement of the integration path in the usual manner to the right over the poles at $\tau=k$, $0\leq k\leq N$, and $\tau=s+k$, $0\leq k\leq M=N-\lfloor\Re (s)\rfloor$, then yields
\[I= \frac{\pi}{\sin \pi s}\left(\frac{\g(a)\g(b)}{\g(n+a+b)}\sum_{k=0}^N \frac{(a)_k(b)_k}{(n+a+b)_k k!}-
\frac{\g(c-a)\g(c-b)}{\g(n+c)\g(1+s)} \sum_{k=0}^M \frac{(c-a)_k(c-b)_k}{(n+c)_k (1+s)_k}\right)+I_N,
\]
where
\[I_N=\frac{1}{2\pi i}\int_{-\infty i+N+\delta}^{\infty i+N+\delta}\frac{\g(a+\tau)\g(b+\tau)}{\g(n+a+b+\tau)\g(1+\tau)}\,
\frac{\pi^2 d\tau}{\sin \pi \tau\,\sin \pi(s-\tau)}\]
with $0<\delta<1$ chosen such that the integration path is a straight line that does not pass through a pole.
From the well-known fact that 
\[\frac{\g(z+a)}{\g(z+b)}\sim z^{a-b} \qquad (|z|\ra\infty,\  |\arg\,z|<\pi),\]
it follows that the ratio of gamma functions appearing in the integrand of $I_N$ is $O(|\tau|^{-n-1})$ as $|\tau|\ra\infty$ on the displaced integration path. It is then readily seen that $I_N=O(N^{-n-1})$ and hence that $I_N\ra 0$ as $N\ra\infty$; the upper limits of summation in the above two series may therefore be replaced by $\infty$. Using the Gauss summation formula (\ref{eGauss}) and expressing the second sum as a ${}_3F_2(1)$ series, we then obtain the expansion in (\ref{e25}). This result can also be obtained directly from the relation between two ${}_3F_2(1)$ series given \cite[Eq.~(4.3.4.3)]{S}.

When $s=0$ ($c=a+b$), all the poles situated at $\tau=k$ are double poles. Putting $\tau=k+\epsilon$, we find that the behaviour of the integrand in (\ref{e28}) as $\epsilon\ra 0$ is
\[\frac{1}{\epsilon^2}\,\frac{\g(a+k)\g(b+k)}{\g(n+a+b+k) k!}\{1+\epsilon[\psi(a+k)+\psi(b+k)-\psi(n+a+b+k)-\psi(1+k)]+\cdots\},\]
with the residue
\[-\frac{\g(a+k)\g(b+k)}{\g(n+a+b+k) k!}\{\psi(n+a+b+k)+\psi(1+k)-\psi(a+k)-\psi(b+k)\}.\]
Then displacement of the integration path over the infinite set of double poles yields the expansion in (\ref{e26}).\ \ \ \ \ $\Box$
\vspace{0.2cm}

Let $m$ be a positive integer. The case $s=m$ requires separate treatment since it is no longer possible to
separate the poles at $\tau=0, 1, \ldots , m-1$ in the integrand in (\ref{e28}), whereas the case $s=-m$ involves both simple and double poles in the integrand.
\begin{theorem}
$\!\!\!$.\ \ Let $m$ and $n$ be  positive integers and $s=c-a-b$, with the quantities $\omega_n$ and $\lambda_n$ as defined in (\ref{e22}). Then, when $s=m$, we have the {\it finite} inverse factorial series
\bee\label{e24}
S_n(a,b;c)=\lambda_n\,\frac{\g(c)\g(c-a-b)}{\g(c-a)\g(c-b)} \sum_{k=0}^{m-1}\frac{(a)_k(b)_k}{(n+a+b)_k k!}
\ee
and when $s=-m$ and either $a$ or $b\neq 1, 2, \ldots , m$ we have
\[S_n(a,b;c)=\frac{\omega_n\g(c)}{m\g(a)\g(b)} \sum_{k=0}^{m-1}\frac{(c-a)_k (c-b)_k}{(n+c)_k (1-m)_k}
+\frac{(-)^m }{m!}\,\frac{\lambda_n \g(c)}{\g(c-a) \g(c-b)}\hspace{2cm}\]
\bee\label{e27}
\hspace{1cm}\times\sum_{k=0}^\infty\frac{(a)_k (b)_k}{(n+a+b)_k k!}
\{\psi(n+a+b+k)+\psi(1+k)-\psi(a+k)-\psi(b+k)\}.
\ee
\end{theorem}
\vspace{0.2cm}

\noindent{\it Proof\,.}\ From (\ref{e23}), we obtain when $s=c-a-b=m$
\[S_n(a,b;c)=\frac{\g(n+a)\g(n+b)}{\g(n)\g(n+a+b)}\,{}_3F_2\left[\begin{array}{c}a, b, f+m-1\\c, f\end{array};1\right],\qquad f:=n+a+b.\]
The above ${}_3F_2(1)$ series can be evaluated by the generalised Karlsson-Minton summation theorem \cite{MS}, \cite[Thm. 6]{MP1} to find
\[{}_3F_2\left[\begin{array}{c}a, b, f+m-1\\c, f\end{array};1\right]=\frac{\g(c)\g(c-a-b)}{\g(c-a)\g(c-b)} \sum_{k=0}^{m-1} (-)^k{\cal A}_k\,\frac{(a)_k (b)_k}{(1-m)_k},\]
where the coefficients ${\cal A}_k$ are given by \cite{MP2} 
\[{\cal A}_k=\frac{(-)^k}{k!}\,{}_2F_1\left[\begin{array}{c}-k, f+m-1\\f\end{array};1\right]=\frac{(-)^k (1-m)_k}{k! (f)_k}.\]
The ${}_2F_1(1)$ evaluation in the above expression follows from Vandermonde's theorem \cite[p.~243]{S}.
Hence we obtain the evaluation when $s=m$ given in (\ref{e24}).
We note that $\lambda_n=1+O(n^{-1})$ as $n\ra\infty$ and the limiting value of the right-hand side of (\ref{e24}) then correctly reduces to the Gauss summation formula in (\ref{eGauss}). The case $m=1$ is seen to agree with
the expression given in \cite[p.~84]{S} when $n$ is replaced by $n+1$.

Finally, when $s=-m$, and provided $a$, $b\neq 1, 2, \ldots ,m$, there are simple poles at $\tau=-1, -2, \ldots , -m$ and double poles at $\tau=k$, $k\geq 0$ situated on the right of the indented integration path in (\ref{e28}). Straightforward evaluation of the residues as in Theorem 1 then yields the expansion in (\ref{e27}).
\ \ \ \ \ $\Box$

\vspace{0.2cm}

\noindent{\bf Remark.}\ \ \ We observe that in the case $s=-m$ when either $a$ or $b=1, 2, \ldots , m$, then
$c-b$ or $c-a=0, -1, -2, \ldots , 1-m$, respectively, and the second sum in (\ref{e27}) accordingly vanishes. Although the integral representation (\ref{e28}) fails in these cases (since the integration path cannot be made to separate the sequences of poles), we conjecture from (\ref{e27}) when $s=-m$ and $a$ or $b=p$, where $p=1, 2, \ldots , m$,
that
\bee\label{e29} 
S_n(a,b;c)=\frac{\omega_n\g(c)}{m\g(a)\g(b)} \sum_{k=0}^{m-p}\frac{(a-m)_k (b-m)_k}{(n+c)_k (1-m)_k}\qquad(s=-m,\ \ p=1, 2, \ldots , m).
\ee
This assertion is supported by numerical evidence.
\vspace{0.2cm}

The results in (\ref{e25}), (\ref{e26}) and (\ref{e27}) involve absolutely convergent series of inverse factorial type in the summation index $n$. This makes these formulas suitable for calculation when $n$ is large.
Using the fact that $\omega_n\sim n^{a+b-c}$ as $n\ra\infty$, we find from (\ref{e25}), (\ref{e24}) and (\ref{e27}) the leading large-$n$ behaviour given by, 
\bee\label{e213a}
S_n(a,b;c)\sim\left\{\begin{array}{ll} \dfrac{\g(c)\g(c-a-b)}{\g(c-a)\g(c-b)} & \Re (s)>0\\
\dfrac{\g(c)\, n^{a+b-c}}{(-\sigma)\g(a)\g(b)} & \Re (s)<0\end{array}\right.
\ee
and, when $s=0$, we have from (\ref{e26}) and the fact that $\psi(n+a+b)\sim \log\,n$  from (\ref{e13a})
\bee\label{e213b}
S_n(a,b;a+b)\sim\frac{\g(a+b)}{\g(a)\g(b)}\,\log\,n.
\ee
The result (\ref{e213a}) when $\Re (s)>0$ corresponds to the case where $S_n(a,b;c)$ converges to a finite sum as $n\ra\infty$ given by the well-known Gauss summation formula (\ref{eGauss}). The leading behaviour when $\Re (s)<0$ and $s=0$, where the sum $S_n(a,b;c)$ diverges as $n\ra\infty$, agrees with that obtained by Hill  \cite{H1, H2}, who derived only the leading terms in the expansions (\ref{e25}), (\ref{e26}) and (\ref{e27}) by means of elaborate algebraic manipulation and induction arguments.

With the help of (\ref{eGauss}) and the properties of the $\psi$-function, the case $s=0$ ($c=a+b$) in (\ref{e26}) can be written in the alternative form
\[S_n(a,b; a+b)=\frac{\g(a+b)}{\g(a)\g(b)}\, \psi(n+a+b)+c_0(a,b)\hspace{4.5cm}\]
\bee\label{e212}
\hspace{3cm}+\frac{\lambda_n\g(a+b)}{\g(a)\g(b)}\sum_{k=1}^\infty \frac{(a)_k (b)_k}{(n+a+b)_k k!}\left\{\sum_{r=0}^{k-1}\frac{1}{n+a+b+r}-\sigma_k(a,b)\right\},
\ee
where
\bee\label{e210}
c_0(a,b)=\frac{\g(a+b)}{\g(a)\g(b)}\{\psi(1)-\psi(a)-\psi(b)\}
\ee
and the coefficients
\bee\label{e211}
\sigma_k(a,b)=\sum_{r=0}^{k-1}\left(\frac{1}{a+r}+\frac{1}{b+r}-\frac{1}{r+1}\right).
\ee

\vspace{0.6cm}

\begin{center}
{\bf 3. \  The Landau constants}
\end{center}
\setcounter{section}{3}
\setcounter{equation}{0}
\renewcommand{\theequation}{\arabic{section}.\arabic{equation}}
From (\ref{e11}), it is seen that the Landau constants $G_n$ are given by
\[G_{n-1}=S_{n}(\fs, \fs ; 1);\]
this corresponds  to the logarithmic case in Theorem 1 with the parametric excess $s=0$. From (\ref{e26}) and (\ref{e212}), we therefore obtain
\[S_n(\fs, \fs ;1)\equiv S_n=\frac{\lambda_n}{\pi} \sum_{k=0}^\infty \frac{(\fs)_k (\fs)_k}{(n+1)_k k!}\{\psi(n+1+k)+\psi(1+k)-2\psi(\fs+k)\}\]
\bee\label{e31}
=\frac{1}{\pi}\psi(n+1)+c_0+\frac{\lambda_n}{\pi}\sum_{k=1}^\infty\frac{(\fs)_k(\fs)_k}{(n+1)_k k!}\left\{\sum_{r=1}^k\frac{1}{n+r}-\sigma_k\right\},
\ee
where
\[\lambda_n=\frac{\g^2(n+\fs)}{\g(n)\g(n+1)},\qquad \sigma_k\equiv\sigma_k(\fs, \fs)=\sum_{r=1}^k\frac{2r+1}{2r-1}\cdot \frac{1}{r}\]
and
\[c_0\equiv c_0(\fs, \fs)=\frac{1}{\pi}\{\psi(1)-2\psi(\fs)\}=\frac{1}{\pi}(\gamma+4\log\,2)\]
with $\gamma$ being the Euler-Mascheroni constant. The first expansion in (\ref{e31}) was given in an equivalent form by Watson \cite[\S 3]{W}, who obtained it by expressing $G_n$ as an integral over $[0,1]$ involving the complete elliptic integral.

Let $M$ be a fixed positive integer. Then, since $\sum_{r=1}^k (n+r)^{-1}<k/(n+1)$ and $\sigma_k<3k$, the remainder $R_M$ after $M-1$ terms in the sum appearing in (\ref{e31}) satisfies the bound\footnote{Watson \cite{W} obtained $|R_M|=O(n^{1-M})$ but this resulted from his use of the crude bound $\psi(x)<x$ for $x>0$.}
\begin{eqnarray}
|R_M|&=&\frac{\lambda_n}{\pi}\sum_{k=M}^\infty \frac{(\fs)_k(\fs)_k}{(n+1)_k k!}\left|\left\{\sum_{r=1}^k\frac{1}{n+r}-\sigma_k\right\}\right|\leq\frac{\lambda_n}{\pi}(3+(n+1)^{-1})
\sum_{k=M}^\infty\frac{(\fs)_k(\fs)_k}{(n+1)_k (k-1)!}\nonumber\\
&\leq&\frac{4\lambda_n}{\pi}\sum_{k=0}^\infty\frac{(\fs)_{k+M}(\fs)_{k+M}}{(n+1)_{k+M}(k+M-1)!}
\leq\frac{4\lambda_n}{\pi}\,\frac{(\fs)_M (\fs)_M}{(n+1)_M}\sum_{k=0}^\infty\frac{(\fs+M)_k (\fs+M)_k}{(n+M+1)_k k!}\nonumber\\
&=&\frac{4\lambda_n}{\pi^2}\,\frac{\g(n)\g^2(M+\fs)}{\g(n-M)}=O(n^{-M})\label{e32}
\end{eqnarray}
as $n\ra\infty$, where we have made use of the facts that 
\[(a)_{k+M}=(a+M)_k (a)_M\]
and $\lambda_n=1+O(n^{-1})$. Hence we have
\bee\label{e33}
S_n=\frac{1}{\pi}\psi(n+1)+c_0+\frac{\lambda_n}{\pi}\sum_{k=1}^{M-1}\frac{(\fs)_k(\fs)_k}{(n+1)_k k!}\left\{\sum_{r=1}^k\frac{1}{n+r}-\sigma_k\right\}+O(n^{-M}).
\ee
\vspace{0.3cm}

\noindent 3.1\ \ {\it Alternative expression for the coefficients appearing in (\ref{e33})}
\vspace{0.3cm}

\noindent We consider the double sum appearing in (\ref{e33}), namely
\[F\equiv \frac{\lambda_n}{\pi}\sum_{k=1}^{M-1}\frac{(\fs)_k(\fs)_k}{(n+1)_k k!}\sum_{p=1}^k\frac{1}{n+p},\]
which we shall rearrange into a single sum.
We express the quantity $(n+p)^{-1}$, where $p$ denotes a positive integer, as an inverse factorial series in the form
\[\frac{1}{n+p}=\sum_{r=1}^p\frac{(-)^{r-1}}{(n+1)_r}\,\prod_{m=1}^{r-1}(p-m),\]
where an empty product is interpreted as unity; compare \cite[p.~177]{PK}. Then we have
\[\sum_{p=1}^k\frac{1}{n+p}=\sum_{p=1}^k \sum_{r=1}^p \frac{(-)^{r-1}}{(n+1)_r}\,\prod_{m=1}^{r-1}(p-m)
=\sum_{r=1}^k \frac{(-)^{r-1}}{(n+1)_r}\sum_{p=r}^k \prod_{m-1}^{r-1}(p-m)\]
\[=\sum_{r=1}^k\frac{(-)^{r-1}}{r(n+1)_r}\,k(k-1)\ldots (k-r+1),\]
since
\[\sum_{p=r}^k \prod_{m=1}^{r-1}(p-m)=\frac{1}{r} \,k(k-1)\ldots (k-r+1).\]
Substitution into $F$ then yields
\[F=\frac{\lambda_n}{\pi}\sum_{k=1}^{M-1}\frac{(\fs)_k (\fs)_k}{(n+1)_k} \sum_{r=1}^k\frac{(-)^{r-1}}{r(n+1)_r (k-r)!}
=\frac{\lambda_n}{\pi}\sum_{r=1}^{M-1}\frac{(-)^{r-1}}{r(n+1)_r}\sum_{k=r}^{M-1}\frac{(\fs)_k (\fs)_k}{(n+1)_k (k-r)!}\]
upon reversal of the order of summation. The upper limit in the inner sum may be replaced by $\infty$, which by the same argument employed to obtain (\ref{e32}) is easily seen to introduce an error term of $O(n^{-M})$ as $n\ra\infty$ with $M$ fixed. Then we find
\begin{eqnarray}
F&=&\frac{\lambda_n}{\pi}\sum_{r=1}^{M-1}\frac{(-)^{r-1}}{r(n+1)_r}\sum_{k=r}^\infty\frac{(\fs)_k(\fs)_k}{(n+1)_k(k-r)!}+O(n^{-M-1})\nonumber\\
&=&\frac{\lambda_n}{\pi}\sum_{r=1}^{M-1}\frac{(-)^{r-1}}{r(n+1)_r}\,\frac{(\fs)_r(\fs)_r}{(n+1)_r} \sum_{k=0}^\infty\frac{(\fs+r)_k(\fs+r)_k}{(n+r+1)_k k!}+O(n^{-M-1})\nonumber\\
&=&\frac{1}{\pi}\sum_{r=1}^{M-1}\frac{(-)^{r-1}}{r(n+1)_r}\,\frac{(\fs)_r(\fs)_r \g(n-r)}{\g(n)}+O(n^{-M-1})\nonumber\\
&=&\frac{1}{\pi}\sum_{r=1}^{M-1}\frac{(-)^{r-1} (\fs)_r (\fs)_r}{r(n^2-1^2)\ldots (n^2-r^2)}+O(n^{-M-1})\label{e34}
\end{eqnarray}
upon use of the Gauss summation formula in (\ref{eGauss}) to evaluate the inner infinite sum.

Combining (\ref{e34}) with (\ref{e33}), we then obtain the following result.
\begin{theorem}$\!\!\!$.\ \ 
Let $M$ be a fixed positive integer. Then, with $S_n\equiv S_n(\fs, \fs; 1)$, we have as $n\ra\infty$
the expansion
\bee\label{e35}
S_n=\frac{1}{\pi}\psi(n+1)+c_0+\frac{1}{\pi}\sum_{r=1}^{K-1}\frac{(-)^{r-1} (\fs)_r (\fs)_r}{r(n^2-1^2)\ldots (n^2-r^2)}-\frac{\lambda_n}{\pi}\sum_{k=1}^{M-1}\frac{(\fs)_k(\fs)_k\sigma_k}{(n+1)_k k!}+O(n^{-M}),
\ee
where $K=\lfloor \fs(M+1)\rfloor$ and the coefficients $\sigma_k$ satisfy the recurrence
\[\sigma_k=\sigma_{k-1}+\frac{2k+1}{2k-1}\cdot\frac{1}{k}\qquad (k\geq 2),\ \ \ \ \sigma_1=3.\]
\end{theorem}
\vspace{0.3cm}

\noindent 3.2\ \ {\it Asymptotic expansion for $S_n$ as $n\ra\infty$}
\vspace{0.3cm}

\noindent Although the expansion (\ref{e35}) is suitable for computation when $n$ is large,
we can obtain the asymptotic expansion of $S_n$ in inverse powers of $n$ by routine algebra. By making use of the facts that
\[\sigma_1=3,\quad \sigma_2=\f{23}{6},\quad\sigma_3=\f{43}{10},\quad\sigma_4=\f{647}{140},\quad\sigma_5=\f{6131}{1260},\quad\sigma_6=\f{70171}{13860},\ldots    \]
and from the expansion of the ratio of two gamma functions \cite[p.~141]{DLMF} that
\[\lambda_n=1-\frac{1}{4n}+\frac{1}{32n^2}+\frac{1}{128n^3}-\frac{5}{2048n^4}-\frac{23}{8192n^5}+O(n^{-6})\qquad (n\ra\infty),\]
we find with the help of {\it Mathematica} that
%\[S_n\sim\frac{1}{\pi}\psi(n+1)+c_0\hspace{9cm}\]
\bee\label{e36}
S_n\sim\frac{1}{\pi}\psi(n+1)+c_0+\frac{1}{\pi}\sum_{k=1}^\infty \frac{(-)^kC_k}{n^k}
\ee
as $n\ra\infty$, where
\[C_1=\f{3}{4},\quad C_2=\f{7}{64},\quad C_3=-\f{3}{128},\quad C_4=-\f{91}{8192},\quad C_5=\f{75}{8192},\quad
C_6=\f{641}{131072}, \ldots \ .\]

If $n$ is replaced by $n+1$ in (\ref{e36}) and use made of the result $\psi(n+2)\sim \log\,(n+1)+\fs (n+1)^{-1}+\f{1}{12}(n+1)^{-2}+ \cdots\ $, then Watson's expansion in (\ref{e13}) is recovered. Similarly, if we
put $h=1$ and replace $n$ by $n-1$ in (\ref{e15}) and make use of (\ref{e13a}), we find agreement with the expansion obtained by Nemes.

\vspace{0.6cm}

\begin{center}
{\bf 4. \  The general case when $s=0$ or a negative integer}
\end{center}
\setcounter{section}{4}
\setcounter{equation}{0}
\renewcommand{\theequation}{\arabic{section}.\arabic{equation}}
The same procedure can be brought to bear on the general logarithmic case $s=0$. From (\ref{e212}),
we have
\[S_n(a,b; a+b)=\frac{\g(a+b)}{\g(a)\g(b)}\, \psi(n+a+b)+c_0(a,b)\hspace{4.5cm}\]
\bee\label{e41}
\hspace{3cm}+\frac{\lambda_n\g(a+b)}{\g(a)\g(b)}\sum_{k=1}^\infty \frac{(a)_k (b)_k}{(n+a+b)_k k!}\left\{\sum_{r=0}^{k-1}\frac{1}{n+a+b+r}-\sigma_k(a,b)\right\},
\ee
where 
\[\lambda_n=\frac{\g(n+a)\g(n+b)}{\g(n)\g(n+a+b)}=1-\frac{ab}{n}+\frac{ab}{2n^2}(a+b-1+ab)+O(n^{-3})\]
as $n\ra\infty$ and $c_0(a,b)$ and the coefficients $\sigma_k(a,b)$ are defined in (\ref{e210}) and (\ref{e211}).
The double sum appearing in (\ref{e41}) can be rearranged, if so desired, following the method used in Section 3 to find
\bee\label{e42}
{\cal F}\equiv \sum_{k=1}^\infty \frac{(a)_k (b)_k}{(n+a+b)_k k!}\sum_{r=0}^{k-1}\frac{1}{n+a+b+r}\hspace{5cm}\]
\[\hspace{4cm}
=\sum_{r=1}^{K-1}\frac{(-)^{r-1} (a)_r(b)_r}{r(n+a+b)_r (n-1)\ldots (n-r)}+O(n^{-M-1}),
\ee
where $K=\lfloor\fs(M+1)\rfloor$. 

The expansion in inverse powers of $n$ then follows by observing that\footnote{This expansion of ${\cal F}$ can be obtained from either side of (\ref{e42}).}
\[{\cal F}=\frac{ab}{n^2}-\frac{a+b-1}{n^3}+O(n^{-4})\]
and 
\[\lambda_n\sum_{k=1}^\infty\frac{(a)_k(b)_k \sigma_k(a,b)}{(n+a+b)_k k!}=\frac{1}{n}(a+b-ab)+\frac{1}{4n^2}(a-1)(b-1)(2a+2b+ab)\]
\[-\frac{1}{36n^3}(a-1)(b-1)\{6(2a^2+2b^2-a-b)+ab(8a+8b+2ab+5))\}+O(n^{-4})\]
Then we obtain the following expansion in the logarithmic case $s=0$
\bee\label{e43}
S_n(a,b;a+b)\sim\frac{\g(a+b)}{\g(a)\g(b)} \psi(n+a+b)+c_0(a,b)+\frac{\g(a+b)}{\g(a)\g(b)}
\sum_{k=1}^\infty \frac{(-)^{k-1}A_k}{n^k}
\ee
as $n\ra\infty$, where
\[A_1=ab-a-b,\qquad A_2=\frac{1}{4}\bl\{(a-1)(b-1)(2a+2b+ab)-4ab\br\},\]
\[A_3=\frac{1}{36}\bl\{(a-1)(b-1)\{6(2a^2+2b^2-a-b)+ab(8a+8b+2ab+5)\}-36ab(a+b-1)\br\}.\]
The above expansion in the case $a=b=\fs$, $c=1$ is seen to correctly reduce to the first three terms in (\ref{e36}).

In the case $s=-m$, where $m$ is a positive integer and either $a$ or $b\neq 1, 2, \ldots , m$, we find after a similar rearrangement of (\ref{e27}) that
\[S_n(a,b;c)\sim\frac{\omega_n\g(c)}{m\g(a)\g(b)} \sum_{k=0}^{m-1}\frac{(c-a)_k (c-b)_k}{(n+c)_k (1-m)_k}
+\frac{(-)^m }{m!}\,\frac{\g(c)}{\g(c-a) \g(c-b)}\bl\{\psi(n+a+b) \]
\bee\label{e44}
+\frac{\g(a)\g(b)}{\g(c-a)\g(c-b)}\,c_0(a,b)+\sum_{k=1}^\infty\frac{(-)^{k-1}A_k}{n^k}\br\}
\ee
as $n\ra\infty$.

To conclude, we present the results of numerical calculations to illustrate the accuracy of the expansions in (\ref{e43}) and (\ref{e44}). In Table 1 below we show the absolute error in the computation of $S_n(a,b;c)$ using the expansion (\ref{e43}) truncated after $k$ terms for selected values of $a$, $b$ and $n$.
\begin{table}[th]
\caption{\footnotesize{Values of the absolute error in the computation of $S_n(a,b;c)$ by 
(\ref{e43}) and (\ref{e44}) for different truncation index $k$. The upper half of the table corresponds to $s=0$ ($c=a+b$) and the lower half to $s=-m$.}}
\begin{center}
\begin{tabular}{r|c|c|c}
\hline
&&&\\[-0.25cm]
\mcol{1}{c|}{$k$} & \mcol{1}{c|}{$a=\f{1}{3},\ b=\f{2}{3}$} & \mcol{1}{c|}{$a=\f{3}{2},\ b=\fs$} & \mcol{1}{c}{$a=\fs+i,\ b=\f{1}{4}$}
\\
&&&\\[-0.35cm]
\mcol{1}{c|}{} & \mcol{1}{c|}{$n=40$} & \mcol{1}{c|}{$n=50$}
& \mcol{1}{c}{$n=100$}
\\
[.15cm]\hline
&&&\\[-0.25cm]
1  & $1.711\times 10^{-5}$      & $2.616\times 10^{-4}$ & $1.545\times 10^{-5}$ \\
2  & $9.618\times 10^{-8}$      & $4.954\times 10^{-6}$ & $1.291\times 10^{-7}$ \\
3  & $\ \,9.845\times 10^{-10}$ & $9.922\times 10^{-8}$ & $1.227\times 10^{-9}$ \\
[.15cm]\hline\hline
&&&\\[-0.25cm]
\mcol{1}{c|}{$k$} & \mcol{1}{c|}{$a=\f{4}{3},\ b=\f{1}{3}$} & \mcol{1}{c|}{$a=\f{3}{2},\ b=-\f{1}{4}$} & \mcol{1}{c}{$a=\f{3}{4}+i,\ b=\f{1}{4}+i$}
\\
&&&\\[-0.35cm]
\mcol{1}{c|}{} & \mcol{1}{c|}{$c=-\f{7}{3},\ n=40$} & \mcol{1}{c|}{$c=\f{1}{4},\ n=50$}
& \mcol{1}{c}{$c=-2+2i,\ n=100$}
\\
[.15cm]\hline
&&&\\[-0.25cm]
1  & $9.820\times 10^{-5}$      & $9.654\times 10^{-6}$      & $6.556\times 10^{-5}$ \\
2  & $1.601\times 10^{-6}$      & $6.888\times 10^{-7}$      & $9.752\times 10^{-7}$ \\
3  & $2.812\times 10^{-8}$      & $\ \,4.141\times 10^{-11}$ & $1.520\times 10^{-8}$ \\
[.15cm]\hline
\end{tabular}
\end{center}
\end{table}

\end{document}